\documentclass[a4paper]{amsart}
\usepackage{amsmath,amssymb}
\usepackage[dvips]{graphics}
\usepackage[all]{xy}

\newtheorem{Theorem}{Theorem}[section]

\newtheorem{Proposition}[Theorem]{Proposition}

\newtheorem{Example}[Theorem]{Example}

\def\labelenumi{\rm ({\makebox[0.8em][c]{\theenumi}})}
\def\theenumi{\arabic{enumi}}
\def\labelenumii{\rm {\makebox[0.8em][c]{(\theenumii)}}}
\def\theenumii{\roman{enumii}}

\renewcommand{\thefigure}
{\arabic{section}.\arabic{figure}}
\makeatletter
\@addtoreset{figure}{section}
\def\@thmcountersep{-}
\makeatother

\numberwithin{equation}{section}

\begin{document}

\title{Knot diagrams on a punctured sphere as a model of string figures}
\author{Masafumi Arai}
\address{Graduate School of Fundamental Science and Engineering, Waseda University, Okubo 3-4-1, Shinjuku, Tokyo 169-8555, Japan (graduated)}
\email{araimasa23@gmail.com}

\author{Kouki Taniyama}
\address{Department of Mathematics, School of Education, Waseda University, Nishi-Waseda 1-6-1, Shinjuku-ku, Tokyo, 169-8050, Japan}
\email{taniyama@waseda.jp}

\thanks{The second author was partially supported by Grant-in-Aid for Scientific Research(A) (No. 16H02145) , Japan Society for the Promotion of Science.}

\subjclass[2020]{Primary 57K10; Secondly 57K20.}

\date{}

\dedicatory{}

\keywords{knot, knot diagram, crossing number, Reidemeister move, string figure, Turaev cobracket, minimal self-intersection number}

\begin{abstract}

A string figure is topologically a trivial knot lying on an imaginary plane orthogonal to the fingers with some crossings. 
The fingers prevent cancellation of these crossings. 
As a mathematical model of string figure we consider a knot diagram on the $xy$-plane in $xyz$-space missing some straight lines parallel to the $z$-axis. These straight lines correspond to fingers. We study minimal number of crossings of these knot diagrams under Reidemeister moves missing these lines. 

\end{abstract}

\maketitle

\section{Introduction}\label{introduction} 

\begin{figure}[htbp]
      \begin{center}
\scalebox{0.5}{\includegraphics*{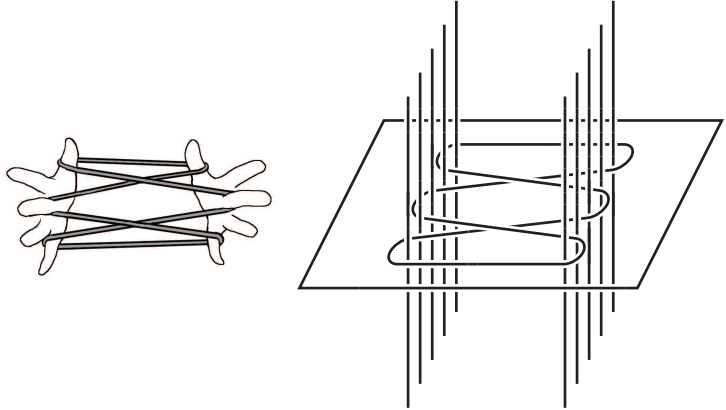}}
      \end{center}
   \caption{A string figure and its mathematical model}
  \label{string-figure}
\end{figure}

A string figure is topologically a trivial knot lying on an imaginary plane orthogonal to the fingers with some crossings. 
The fingers prevent cancellation of these crossings. 
As a mathematical model of string figure we consider a knot lying on ${\mathbb R}^{2}\times\{0\}\setminus{\mathcal L}$ where ${\mathbb R}^{2}\times\{0\}$ is the $xy$-plane in the $xyz$-space ${\mathbb R}^{3}$ and 
${\mathcal L}$ is a union of finitely many straight lines each of which is parallel to the $z$-axis. We identify a knot lying on a plane with a knot diagram on the plane. 
See Figure \ref{string-figure}. 

By the one-point compactification of the pair ${\mathbb R}^{2}\times\{0\}\subset{\mathbb R}^{3}$ we consider the pair ${\mathbb S}^{2}\subset{\mathbb S}^{3}$. Namely we consider a knot $K$ in the $3$-sphere ${\mathbb S}^{3}$ and a diagram $D$ of $K$ on the $2$-sphere ${\mathbb S}^{2}$. Instead of straight lines parallel to the $z$-axis that restrict the deformation of $K$ we remove some open disks from ${\mathbb S}^{2}$ and restrict the deformation of $D$ to Reidemeister moves performed on this punctured sphere. 

The precise formulation is as follows. 
Let $K$ be a knot in ${\mathbb S}^{3}$ and $D$ a diagram of $K$ on ${\mathbb S}^{2}$. 
Let $P(D)\subset {\mathbb S}^{2}$ be the immersed circle obtained from $D$ by forgetting the over/under crossing information at each crossing point of $D$. We denote $P(D)$ by $P$ when the choice of $D$ is clear. 
Sometimes $P$ does not come from $D$ and it is simply an image of a generic immersion $\varphi:{\mathbb S}^{1}\to{\mathbb S}^{2}$. Namely $P=\varphi({\mathbb S}^{1})$. 
Let $N(P)\subset {\mathbb S}^{2}$ be the regular neighbourhood of $P$ in ${\mathbb S}^{2}$. 
Let ${\mathcal R}(P)$ be the set of connected components of ${\mathbb S}^{2}\setminus N(P)$. 
We denote ${\mathcal R}(P)$ by ${\mathcal R}$ when the choice of $P$ is clear. 
Note that each element of ${\mathcal R}$ is an open disk whose closure is a closed disk. Let $S$ be a subset of ${\mathcal R}$. Set $F(S)={\mathbb S}^{2}\setminus\bigcup_{R\in S}R$. Note that $D$ is still a diagram of $K$ on this $|S|$-punctured sphere $F(S)$. 
Note also that $F(\emptyset)={\mathbb S}^{2}$, $F({\mathcal R})=N(P)$ and if $S_{1}\subset S_{2}$ then $F(S_{1})\supset F(S_{2})$. 
Let $C(D)$ be the number of crossing points of $D$. Let $c(D)$ be the minimal number of $C(E)$ where $E$ varies over all knot diagrams obtained from $D$ by Reidemeister moves performed on ${\mathbb S}^{2}$. Namely $c(D)=c(K)$ is the minimal crossing number of the knot $K$. Let $c(D,S)$ be the minimal number of $C(E)$ where $E$ varies over all knot diagrams on $F(S)$ obtained from $D$ by Reidemeister moves performed on $F(S)$. 

\begin{figure}[htbp]
      \begin{center}
\scalebox{0.5}{\includegraphics*{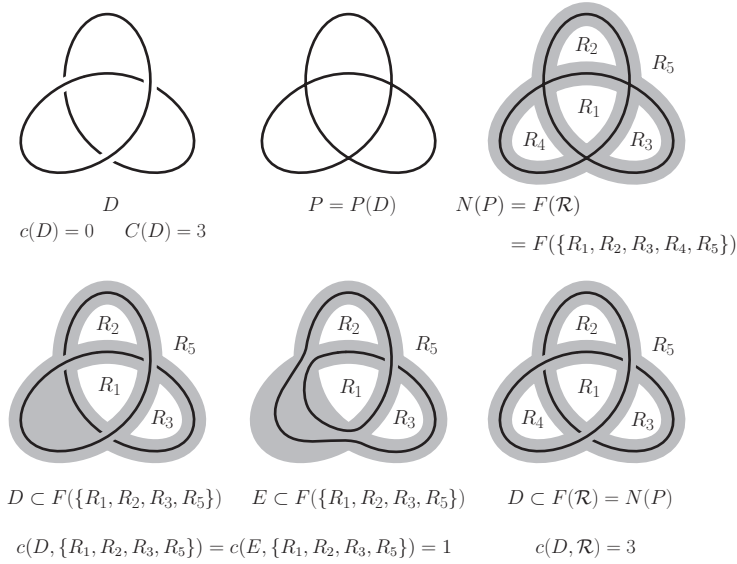}}
      \end{center}
   \caption{An example}
  \label{trivial-diagram}
\end{figure}

\begin{Example}\label{Example}
{\rm
Let $D$ be a knot diagram on ${\mathbb S}^{2}$ of a trivial knot in ${\mathbb S}^{3}$ as illustrated in Figure \ref{trivial-diagram} where ${\mathbb S}^{2}$ is regarded as one-point compactification of ${\mathbb R}^{2}$ and $D$ in ${\mathbb R}^{2}$ is illustrated. Let ${\mathcal R}=\{R_{1},R_{2},R_{3},R_{4},R_{5}\}$ as illustrated in Figure \ref{trivial-diagram}. 
Then we see that $D$ and $E$ in Figure \ref{trivial-diagram} are transformed into each other by second Reidemeister move on $F(\{R_{1},R_{2},R_{3},R_{5}\})$. Therefore $c(D,\{R_{1},R_{2},R_{3},R_{5}\})\leq C(E)=1$. 
It follows from Theorem \ref{theorem4} in Section \ref{proofs} that $c(D,\{R_{1},R_{5}\})\geq1$. 
Since $c(D,\{R_{1},R_{2},R_{3},R_{5}\})\geq c(D,\{R_{1},R_{5}\})$ we have $c(D,\{R_{1},R_{2},R_{3},R_{5}\})\geq1$.
Therefore we have $c(D,\{R_{1},R_{2},R_{3},R_{5}\})=1$. 
By Theorem \ref{theorem} below we have $c(D,{\mathcal R}(P))=C(D)=3$. 
}
\end{Example}

By definition we have the following proposition. 


\begin{Proposition}\label{proposition}

Let $D$ be a knot diagram on ${\mathbb S}^{2}$ and $P=P(D)$. 

{\rm (1)} $c(D,{\emptyset})=c(D)$. 

{\rm (2)} Let $S$ be a subset of ${\mathcal R}(P)$. Then $c(D)\leq c(D,S)\leq C(D)$. 

{\rm (3)} Suppose that $S_{1}\subset S_{2}\subset{\mathcal R}(P)$. Then $c(D,S_{1})\leq c(D,S_{2})$.

\end{Proposition}

Namely the map $S\mapsto c(D,S)$ is an order-preserving map from the power set $2^{{\mathcal R}(P)}$ to the set of all non-negative integers ${\mathbb Z}_{\geq0}$ where $2^{{\mathcal R}(P)}$ is partially ordered by set inclusion. Then it is natural to ask the question whether $c(D,{{\mathcal R}(P)})=C(D)$ or not. We have the following affirmative answer. 

\begin{Theorem}\label{theorem}

Let $D$ be a knot diagram on ${\mathbb S}^{2}$ and $P=P(D)$. 
Then $c(D,{{\mathcal R}(P)})=C(D)$.

\end{Theorem}

We have asked the question above at some opportunities in 2019. An affirmative answer above is first given by K. Tagami. His proof of Theorem \ref{theorem} is based on his result in \cite[Corollary 4.11]{Tagami}. After he told us his proof, we have noticed that there is a simple proof using the Turaev cobracket \cite{Turaev} of Theorem \ref{theorem2} from which Theorem \ref{theorem} immediately follows. 
Then Z. Cheng informed us that Theorem \ref{theorem2} immediately follows from \cite[Theorem 4.2]{H-S}. Since the proof Theorem 4.2 in \cite{H-S} is relatively long, we think that our proof of Theorem \ref{theorem2} based on the Turaev cobracket is worth stating. It is given in Section \ref{proofs}. 

\begin{Theorem}\label{theorem2}

Let $\varphi:{\mathbb S}^{1}\to{\mathbb S}^{2}$ be a generic immersion and $P=\varphi({\mathbb S}^{1})$. 
Let $N(P)$ be a regular neighbourhood of $P$ in ${\mathbb S}^{2}$. Then the minimal self-intersection number of $\varphi$ among all generic immersions homotopic to $\varphi$ on $N(P)$ is equal to the number of crossings of $\varphi$. 

\end{Theorem}

Suppose that $c(D)<C(D)$. Then $c(D)=c(D,{\emptyset})<c(D,{{\mathcal R}(P)})=C(D)$. 
We ask what is the smallest $S$ with $c(D)<c(D,S)$. We also ask what is the largest $S$ with $c(D)=c(D,S)$.
We prepare the following definitions to be more precise. 

Set $m=C(D)$. Then it is well-known that $|{\mathcal R}(P)|=m+2$. Let $n$ be a non-negative integer with $0\leq n\leq m+2$. Define $c_{\rm max}(D,n)$ (resp. $c_{\rm min}(D,n)$) to be the maximum (resp. minimum) of $c(D,S)$ where $S$ varies over all subset of ${\mathcal R}(P)$ with $|S|=n$. 
By definition we have $c_{\rm max}(D,0)=c_{\rm min}(D,0)=c(D,\emptyset)=c(D)$ and $c_{\rm max}(D,m+2)=c_{\rm min}(D,m+2)=c(D,{\mathcal R}(P))=C(D)$. 
Moreover we have the following proposition. 

\begin{Proposition}\label{inequalities}

Let $D$ be a knot diagram on ${\mathbb S}^{2}$ with $C(D)=m$.
Then we have the following inequalities. 

{\rm
\begin{alignat*}{6}
&c_{\rm max}(D,0)&&=c_{\rm max}(D,1)&&\leq c_{\rm max}(D,2)&&\leq&&\cdots\leq c_{\rm max}(D,m+1)&&\leq c_{\rm max}(D,m+2)\\
&\hspace{8mm}\rotatebox{90}{=}&&\hspace{11mm}\rotatebox{90}{=}&&\hspace{10mm}\rotatebox{90}{$\leq$}&&&&\cdots\hspace{12mm}\rotatebox{90}{$\leq$}&&\hspace{13mm}\rotatebox{90}{=}\\
&c_{\rm min}(D,0)&&=c_{\rm min}(D,1)&&=c_{\rm min}(D,2)&&\leq&&\cdots\leq c_{\rm min}(D,m+1)&&\leq c_{\rm min}(D,m+2).
\end{alignat*}
}

\end{Proposition}

A proof of Proposition \ref{inequalities} is given in Section \ref{proofs}.
We pay attention to the inequality $c_{\rm max}(D,0)=c_{\rm max}(D,1)\leq c_{\rm max}(D,2)$ in Proposition \ref{inequalities}. 
Then the next question is whether or not $c_{\rm max}(D,1)= c_{\rm max}(D,2)$. The following is a partial answer to this question. 

\begin{Theorem}\label{theorem3}

Let $D$ be a knot diagram on ${\mathbb S}^{2}$ of a trivial knot with $C(D)>0$. 
Then $c_{\rm max}(D,2)>0$. 

\end{Theorem}

A proof is given in Section \ref{proofs}. 

\vskip 3mm

This paper is based on the graduation thesis of the first author submitted to Waseda University. The basic idea of this paper owes to him. 


\section{Proofs}\label{proofs} 

\vskip 3mm

Let $\varphi:{\mathbb S}^{1}\to{\mathbb S}^{2}$ be a generic immersion and $P=\varphi({\mathbb S}^{1})$. 
Let ${\mathcal C}$ be the set of all crossing points of $P$. Each connected component of $P\setminus{\mathcal C}$ is said to be an {\it edge} of $P$. Let $e$ be an edge of $P$. 
Let $u$ and $v$ be mutually distinct connected components of ${\mathbb S}^{2}\setminus P$ such that each of the closure of $u$ and the closure of $v$ contains $e$. Let $U$ and $V$ be the elements of ${\mathcal R}(P)$ with $U\subset u$ and $V\subset v$. Then we say that $U$ and $V$ are {\it adjacent along $e$}. 
Let $F(e)=F(\{U,V\})={\mathbb S}^{2}\setminus\{U\cup V\}$. Note that $F(e)$ is an annulus. 
Let $S$ be a subset of ${\mathcal R}$. A loop on a surface $F(S)$ is said to be nontrivial on $F(S)$ if it is not null-homotopic on $F(S)$. We give an orientation to ${\mathbb S}^{1}$ and regard $P=\varphi({\mathbb S}^{1})$ as an oriented loop on $N(P)$.

\vskip 3mm

\begin{figure}[htbp]
      \begin{center}
\scalebox{0.5}{\includegraphics*{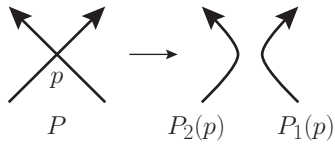}}
      \end{center}
   \caption{Smoothing $P$ at $p$}
  \label{smoothing}
\end{figure}

\noindent{\bf Proof of Theorem \ref{theorem2}.} 
Let $H$ be the set of all free homotopy classes of nontrivial oriented loops on $N(P)$. 
Let $[P]_{F(S)}$ be the free homotopy class of $P$ on $F(S)$ for each $S\subset{\mathcal R}$. For any edge $e$ of $P$ we see that $[P]_{F(e)}$ is primitive. Namely $[P]_{F(e)}$ is not a power of another class. 
Therefore $[P]_{N(P)}$ is also primitive. Therefore $[P]_{N(P)}$ is an element of $H$. 
For each crossing point $p\in{\mathcal C}$ we have two loops $P_{1}(p)$ and $P_{2}(p)$ on $N(P$) as illustrated in Figure \ref{smoothing}. We see that there exist edges $d$ and $e$ of $P$ such that both $[P_{1}(p)]_{F(d)}$ and $[P_{2}(p)]_{F(e)}$ are primitive. Therefore both $[P_{1}(p)]_{N(P)}$ and $[P_{2}(p)]_{N(P)}$ are elements of $H$. 
Let ${\mathbb Z}H$ be the free ${\mathbb Z}$-module generated by $H$. 
Let $\tau:{\mathbb Z}H\to{\mathbb Z}H\otimes_{\mathbb Z}{\mathbb Z}H$ be the Turaev cobracket. 
Then 
\[
\tau([P]_{N(P)})=\sum_{p\in{\mathcal C}}([P_{1}(p)]_{N(P)}\otimes [P_{2}(p)]_{N(P)}-[P_{2}(p)]_{N(P)}\otimes [P_{1}(p)]_{N(P)})
\]
and we see that the minimal self-intersection number of $[P]_{N(P)}$ is greater than or equal to half the number of terms in the linear combination $\tau([P]_{N(P)})$ in ${\mathbb Z}H\otimes_{\mathbb Z}{\mathbb Z}H$. See for example \cite{Cahn}. Suppose that $p\neq q$ or $j\neq k$. Then we see that there is an edge $e$ of $P$ whose closure contains $p$ such that $[P_{j}(p)]_{F(e)}\neq[P_{k}(q)]_{F(e)}$. Therefore $[P_{j}(p)]_{N(P)}\neq[P_{k}(q)]_{N(P)}$. 
Thus we see that $[P_{j}(p)]_{N(P)}=[P_{k}(q)]_{N(P)}$ if and only if $p=q$ and $j=k$. Therefore no terms in the definition of $\tau([P]_{N(P)})$ cancel each other and the number of terms of $\tau([P]_{N(P)})$ is exactly twice the number of elements of ${\mathcal C}$ as desired. 
$\Box$

\noindent{\bf Proof of Proposition \ref{inequalities}.} 
By definition and by Proposition \ref{proposition} (3) we have $c_{\rm max}(D,i)\leq c_{\rm max}(D,i+1)$ and $c_{\rm min}(D,i)\leq c_{\rm min}(D,i+1)$. By definition we have $c_{\rm max}(D,i)\geq c_{\rm min}(D,i)$. We have already remarked that $c_{\rm max}(D,0)=c_{\rm min}(D,0)$ and $c_{\rm max}(D,m+2)=c_{\rm min}(D,m+2)$. 
Therefore it is sufficient to show $c_{\rm max}(D,0)=c_{\rm max}(D,1)$ and $c_{\rm min}(D,0)=c_{\rm min}(D,1)=c_{\rm min}(D,2)$. 

For any $R\in{\mathcal R}$ we see that $F(\{R\})$ is a disk. Note that $D$ is a diagram of $K$ on the disk $F(\{R\})$. 
Therefore $D$ can be deformed into a diagram $E$ of $K$ with $C(E)=c(K)$ by Reidemeister moves on $F(\{R\})$. 
Thus we have $c(D,\{R\})=c(K)=c(D)=c(D,0)$. Therefore $c_{\rm max}(D,0)=c_{\rm max}(D,1)$ as desired. 
Let $P=P(D)$. Let $e$ be any edge of $P$. Let $U$ and $V$ be the elements of ${\mathcal R}(P)$ mutually adjacent along $e$. 
Then $F(e)=F(\{U,V\})$ is an annulus and $D$ is a diagram of $K$ on $F(e)$. Note that $U$ and $V$ may be adjacent not only along $e$ but also along some other edges of $P$. 
Let $k$ be the number of edges of $P$ along which $U$ and $V$ are adjacent. Then we see that $D$ is a diagram-connected sum of $k$ local knot diagrams. 
By an ambient isotopy on $F(e)$ we gather these local knot diagrams in one place. 
Then we see that $D$ can be transformed into a diagram $E$ of $K$ with $C(E)=c(K)$ by Reidemeister moves on $F(e)$. This means that $c(D,\{U,V\})=c(K)=c(D)=c(D,0)$. Thus we have $c_{\rm min}(D,2)=c_{\rm min}(D,0)$. Since $c_{\rm min}(D,0)\leq c_{\rm min}(D,1)\leq c_{\rm min}(D,2)$ we have $c_{\rm min}(D,0)=c_{\rm min}(D,1)=c_{\rm min}(D,2)$ as desired. 
$\Box$

\vskip 3mm

\begin{figure}[htbp]
      \begin{center}
\scalebox{0.5}{\includegraphics*{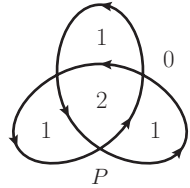}}
      \end{center}
   \caption{An example of Alexander numbering}
  \label{Alexander-numbering}
\end{figure}

Let $\varphi:{\mathbb S}^{1}\to{\mathbb S}^{2}$ be a generic immersion and $P=\varphi({\mathbb S}^{1})$. 
Suppose that $P$ is oriented. A locally constant map $a:{\mathbb S}^{2}\setminus P\to{\mathbb Z}$ is said to be an {\it Alexander numbering} if for each edge $e$ of $P$ the value of $a$ of a point on the left side of $e$ is always one greater than the value of $a$ of a point on the right side of $e$. 
Alexander numbering always exists and unique up to an additive constant. 
For example, for $P\subset{\mathbb R}^{2}\subset{\mathbb R}^{2}\cup\{\infty\}\cong{\mathbb S}^{2}$, the winding number of $P$ around each point in ${\mathbb R}^{2}\setminus P$ defines an Alexander numbering. 
See Figure \ref{Alexander-numbering} for an example. 

\begin{Theorem}\label{theorem4}

Let $\varphi:{\mathbb S}^{1}\to{\mathbb S}^{2}$ be a generic immersion and $P=\varphi({\mathbb S}^{1})$. 
Suppose that $P$ is oriented. Let $a:{\mathbb S}^{2}\setminus P\to{\mathbb Z}$ be an Alexander numbering. 
Let $U$ and $V$ be elements of ${\mathcal R}(P)$, $x\in U$ and $y\in V$. 
Then the minimal self-intersection number of $\varphi$ among all generic immersions homotopic to $\varphi$ on $F(\{U,V\})$ is greater than or equal to $|a(x)-a(y)|-1$. 

\end{Theorem}

\noindent{\bf Proof.} 
Since the case $U=V$ is trivial we may suppose $U\neq V$. Then $F(\{U,V\})$ is an annulus and $P$ is freely homotopic on $F(\{U,V\})$ to a power of the core curve of $F(\{U,V\})$. We see by the definition of Alexander numbering that the power is equal to $a(x)-a(y)$ up to sign. Therefore the curve must have at least $|a(x)-a(y)|-1$ self-intersection points. 
$\Box$

\vskip 3mm

\begin{figure}[htbp]
      \begin{center}
\scalebox{0.5}{\includegraphics*{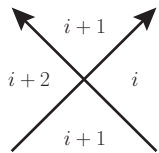}}
      \end{center}
   \caption{Alexander numbering around a crossing point}
  \label{Alexander-numbering2}
\end{figure}

\noindent{\bf Proof of Theorem \ref{theorem3}.}
Suppose that $P=P(D)$ is oriented. 
Let $a:{\mathbb S}^{2}\setminus P\to{\mathbb Z}$ be an Alexander numbering. 
Since $C(D)>0$ there is a crossing point $p$ of $P$. Paying attention to a neighbourhood of $p$ we see that there are elements $U$ and $V$ of ${\mathcal R}(P)$ such that $|a(x)-a(y)|\geq2$ for $x\in U$ and $y\in V$. See Figure \ref{Alexander-numbering2}. Then by Theorem \ref{theorem4} we have $c(D,\{U,V\})\geq1$. Therefore $c_{\rm max}(D,2)\geq c(D,\{U,V\})\geq1$ as desired. 
$\Box$

\section{Future directions}\label{future} 

After Proposition \ref{inequalities} and Theorem \ref{theorem3} we are interested in the first $i$ with $c_{\rm max}(D,i)< c_{\rm max}(D,i+1)$, the first $i$ with $c_{\rm min}(D,i)< c_{\rm min}(D,i+1)$, the last $i$ with $c_{\rm max}(D,i)< c_{\rm max}(D,i+1)$, and the last $i$ with $c_{\rm min}(D,i)< c_{\rm min}(D,i+1)$. 
We define the following to be more precise. 
Let $\alpha(D)=\min\{n\mid c_{\rm max}(D,n)>c(D)\}$, $\beta(D)=\min\{n\mid c_{\rm min}(D,n)>c(D)\}$, $\gamma(D)=\max\{n\mid c_{\rm max}(D,n)<C(D)\}$, and $\delta(D)=\max\{n\mid c_{\rm min}(D,n)<C(D)\}$. The decision of these numbers will be a problem. 
For example, the following theorem is an immediate consequence of Proposition \ref{inequalities} and Theorem \ref{theorem3}. 

\begin{Theorem}\label{theorem3-1}

Let $D$ be a knot diagram on ${\mathbb S}^{2}$ of a trivial knot with $C(D)>0$. 
Then $\alpha(D)=2$. 

\end{Theorem}

\begin{figure}[htbp]
      \begin{center}
\scalebox{0.5}{\includegraphics*{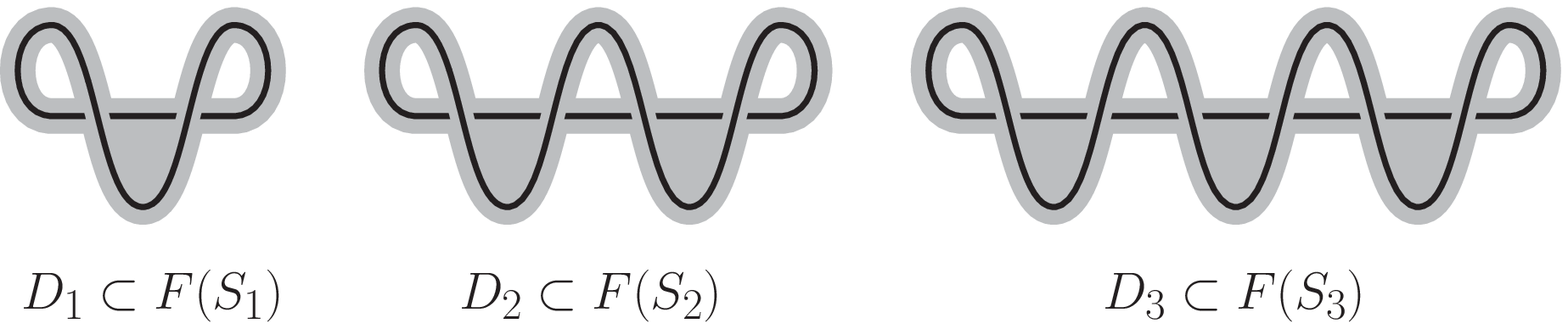}}
      \end{center}
   \caption{Examples}
  \label{trivial-diagram2}
\end{figure}

\begin{Example}\label{example2}
{\rm
Let $D_{n}$ be a knot diagram on a surface $F(S_{n})$ as illustrated in Figure \ref{trivial-diagram2}. 
Note that $C(D_{n})=2n$ and $|S_{n}|=n+2$. By second Reidemeister moves we have $c(D_{n},S_{n})=0$. 
By an argument using Alexander numbering analogous to that of the proof of Theorem \ref{theorem3}, we have $\beta(D_{n})=n+3$. 
Therefore $\displaystyle{\beta(D_{n})=\frac{C(D_{n})}{2}+3}$. 
}
\end{Example}

\vskip 3mm

\section*{Acknowledgments} The authors are grateful to Professor Keiji Tagami for showing us his proof of Theorem \ref{theorem}. 
We thank Professor Zhiyun Cheng for his comments on the proof of Theorem \ref{theorem2}. 
We are also grateful to Ms. Miki Arai for drawing the left figure of Figure \ref{string-figure}. 
We thank the referee for his/her valuable comments.

\vskip 3mm

{\normalsize

\begin{thebibliography}{99}

\bibitem{Cahn}{P. Cahn}, {A generalization of the Turaev cobracket and the minimal self-intersection number of a curve on a surface}, {\it New York J. Math.}, {\bf 19} (2013), 253-283. 

\bibitem{H-S}{J. Hass and P. Scott}, {Intersections of curves on surfaces}, {Israel J. Math.}, {\bf 51} (1985), 90-120.

\bibitem{Tagami}{K. Tagami}, {A Khovanov type invariant derived from an unoriented HQFT for links in thickened surfaces}, {\it Internat. J. Math.}, {\bf 24} (2013), 1350078, 28 pp. 

\bibitem{Turaev}{V. Turaev}, {Skein quantization of Poisson algebras of loops on surfaces}, {\it Ann. Sci. \'{E}cole Norm. Sup. (4)}, {\bf 24} (1991), 635-704. 

\end{thebibliography}
\end{document}